\documentclass[12pt]{amsart}
\usepackage{amsmath,amssymb,amsthm,amscd,graphics,psfrag
}
\newcommand\Ad{\operatorname{Ad}}

\newcommand\C{{\mathbb C}}
\newcommand\caff{{{\mathrm c}_{\mathrm a}}}

\newcommand\calE{{\mathcal E}}
\newcommand\calO{{\mathcal O}}
\newcommand\CB{\operatorname{CB}}
\newcommand\CC{\operatorname{CC}}
\newcommand\eps{\varepsilon}
\newcommand\F{{\mathcal F}}
\newcommand\g{{\mathfrak g}}
\newcommand\gout{\g_{\mathrm{out}}}
\newcommand\goutX{{\g_{\X,\mathrm{out}}}}
\newcommand\gtw{\g^{\mathrm{tw}}}
\newcommand\gtwX{\g^{\mathrm{tw}}_\X}

\newcommand\hatCC{{\widehat\sCC}}
\newcommand\hatM{{\widehat\M}}
\newcommand\hatOS{{\hat\calO_{S/S_0}}}
\newcommand\h{{\mathfrak h}}
\newcommand\HOM{\mathop{\mathcal{H}\mathit{om}}\nolimits}
\newcommand\Hom{\operatorname{Hom}\nolimits}
\newcommand\id{{\mathrm{id}}}
\renewcommand\Im{\operatorname{Im}\nolimits}
\newcommand\Ind{\operatorname{Ind}\nolimits}
\newcommand\isoto{\overset\sim\rightarrow}
\newcommand\isoot{\overset\sim\leftarrow}
\newcommand\K{{\mathcal K}}
\newcommand\khat{{\hat k}} 
\renewcommand\L{{\mathcal L}}
\newcommand\M{{\mathcal M}}
\newcommand\m{{\mathfrak m}}
\newcommand\N{{\mathbb N}}
\newcommand\n{{\mathfrak n}}
\newcommand\orb{\mathrm{orb}}

\newcommand\Proj{{\mathbb P}}
\newcommand\Res{\mathop{\operatorname{Res}}\nolimits}
\newcommand\sCB{\mathcal{CB}}
\newcommand\sCBorb{\mathcal{CB}_\orb}
\newcommand\sCBtri{\mathcal{CB}_\trig}
\newcommand\sCC{\mathcal{CC}}
\newcommand\sCCorb{\mathcal{CC}_\orb}
\newcommand\sCCtri{\mathcal{CC}_\trig}
\renewcommand\setminus{\smallsetminus}
\newcommand\simeqq{\cong}
\newcommand\tensor{\otimes}
\newcommand\tr{\operatorname{tr}}
\newcommand\trig{\mathrm{trig}}

\newcommand\wadd{w^{\mathrm{add}}}
\newcommand\wmul{w^{\mathrm{mul}}}
\newcommand\wt{\operatorname{wt}}
\newcommand\X{{\mathfrak X}}

\newcommand\Z{{\mathbb Z}}
\newtheorem{theorem}{Theorem}[section]
\newtheorem{lemma}[theorem]{Lemma}
\newtheorem{proposition}[theorem]{Proposition}
\newtheorem{corollary}[theorem]{Corollary}
\theoremstyle{definition}
\newtheorem{definition}[theorem]{Definition}
\theoremstyle{remark}

\newtheorem{rem}[theorem]{Remark}
\newcommand\thmref[1]{Theorem~\ref{#1}}
\newcommand\propref[1]{Proposition~\ref{#1}}
\newcommand\secref[1]{\S\ref{#1}}

\newcommand\lemref[1]{Lemma~\ref{#1}}

\newcommand\defref[1]{Definition~\ref{#1}}

\begin{document}
%
\title[Trigonometric Degeneration and Orbifold WZW]
{Trigonometric Degeneration\\ and\\
Orbifold Wess-Zumino-Witten Model. II
}
\author{Takashi Takebe}
\address{%
Department of Mathematics\\
Ochanomizu University\\
Otsuka 2-1-1, Bunkyo-ku\\
Tokyo, 112-8610, Japan
}

\email{takebe@math.ocha.ac.jp}

\thanks{This work is partly supported by the Grant-in-Aid for
Scientific Research (C) of the Japan Society for the Promotion of
Science, No.\ 15540014.}

\dedicatory{Dedicated to Professor Akihiro Tsuchiya on his 60th birthday.} 

\begin{abstract}
 The sheaves of conformal blocks and conformal coinvariants of the
 twisted WZW model have a factorisation property and are locally free
 even at the boundary of the moduli space, where the elliptic KZ
 equations and the Baxter-Belavin elliptic $r$ matrix degenerate to the
 trigonometric KZ equations and the trigonometric $r$  matrix,
 respectively. Etingof's construction of the elliptic KZ equations is
 geometrically interpreted.
\end{abstract}

\maketitle

\section{Introduction}
\label{sec:introduction}

This is a continuation of the paper \cite{tak:04}. We showed there that
the trigonometric WZW model is factorised into the orbifold WZW
models. Using this result, we show in the present article that the
trigonometric WZW model is indeed the degenerate twisted WZW models on
elliptic curves defined in \cite{kur-tak:97}. More precisely, we prove
that there are locally free sheaves over the partially compactified
family of elliptic curves, the fibre of which are the space of conformal
blocks or the space of conformal coinvariants of the twisted WZW model
at a generic point and those of the trigonometric WZW model at the
discriminant locus, when all inserted modules are either Weyl modules
(\propref{prop:coherence}) or integrable highest weight modules
(\thmref{thm:loc-free}).

Since the elliptic $r$ matrix (cf.\ \cite{bel-dri:82}, \cite{eti:94})
describing the elliptic KZ equations degenerates to the trigonometric
$r$ matrix, the above fact is naturally expected, though rigorous
proof requires careful algebro-geometric arguments as in \cite{tuy:89}.

The paper is organised as follows. After reviewing the twisted WZW model
on elliptic curves in \secref{sec:twwzw-ell} to recall basic notions and
notations, we define family of elliptic curves with a singular fibre and
a twisted Lie algebra bundle over it in \secref{sec:family}. In
\secref{sec:cc,cb} main objects of this paper, the sheaves of conformal
coinvariants and conformal blocks, are defined and their coherence is
proved. In particular when all the modules inserted to the curve are
Weyl modules, they are locally free.  To prove the locally freeness of
the sheaf of conformal coinvariants for integrable highest weight
modules, we examine its behaviour at the discriminant locus. In this
case the factorisation theorem, Theorem 7 of \cite{tak:04}, is refined
in \secref{sec:sheaf-tri/orb}. The proof of locally freeness in
\secref{sec:loc-free} follows the strategy of \cite{tuy:89} and
\cite{tsu-kuw:04}.

\subsection*{Notations}

We use the following notations besides other ordinary conventions in
mathematics. 

\begin{itemize}
 \item $N$, $L$: fixed integers. $N \geqq 2$ will be the matrix size and
       $L \geqq 1$ will be the number of the marked points on a curve.

 \item $C_N := \Z/N\Z$: the cyclic group of order $N$.

 \item When $X$ is an algebraic variety, $\calO_X$ denotes the structure
       sheaf of $X$. When $P$ is a point on $X$ and $\F$ is an
       $\calO_X$-sheaf, $\F_P$ denotes the stalk of $\F$ at $P$. $\m_P$ is
       the maximal ideal of the local ring $\calO_{X,P}$.
       $\F|_P:=\F_P/\m_P\F_P$, $\F^\wedge_P = \projlim_{n\to\infty}
       \F_P/\m_P^n \F_P$ are the fibre of $\F$ at $P$ and the
       $\m_P$-adic completion of $\F_P$, respectively.

 \item We shall use the same symbol for a vector bundle and for a
       locally free $\calO_X$-module consisting of its local holomorphic
       sections.

\end{itemize}

\section{Twisted WZW model on elliptic curves}
\label{sec:twwzw-ell}

In this section we briefly review the twisted WZW model on elliptic
curves. See \cite{kur-tak:97} for details.

We fix an invariant inner product of $\g = sl_N(\C)$ by
\begin{equation}
    (A | B) := \tr (AB) \quad \text{for $A, B \in \g$}.
\label{def:inner-product}
\end{equation}
Define matrices $\beta$ and $\gamma$ by
\begin{equation}
    \beta :=  \begin{pmatrix}
         0 & & 1 & \ \ &       & \ \ & 0 \\
           & & 0 & \ \ &\ddots & \ \ &   \\
           & &   & \ \ &\ddots & \ \ & 1 \\
         1 & &   & \ \ &       & \ \ & 0
         \end{pmatrix},
    \quad
    \gamma := \begin{pmatrix}
         1 &            &        & 0 \\
           & \eps^{-1}  &        &   \\
           &            & \ddots &   \\
         0 &            &        & \eps^{1-N}
         \end{pmatrix},
\label{def:beta,gamma}
\end{equation}
where $\eps = \exp (2\pi i/ N)$.  Then we have $\beta^N=\gamma^N=1$
and $\gamma\beta=\eps\beta\gamma$.

Let $E = E_\tau$ be the elliptic curve with modulus $\tau$: $E_\tau :=
\C/\Z + \tau \Z$. We define a Lie algebra bundle $\gtw$ with fibre $\g =
sl_N(\C)$ over $E$ by
\begin{equation}
    \gtw := (\C \times \g)/{\approx},
\label{def:gtw}
\end{equation}
where the equivalence relations $\approx$ are defined by
\begin{equation}
   (z, A) \approx (z+1, \Ad\gamma (A)) 
          \approx (z+\tau, \Ad\beta (A)).
\label{def:eq-rel-gtw}
\end{equation}
Let $J_{ab} = \beta^a \gamma^{-b}$, which satisfies
\begin{equation}
    \Ad \gamma (J_{ab}) = \eps^a J_{ab},\qquad 
    \Ad \beta (J_{ab}) = \eps^b J_{ab}.
\label{Ad(*)Jab}
\end{equation}
Global meromorphic sections of $\gtw$ are linear combinations of $J_{ab}
f(z)$ ($a,b = 0,\dots,N-1$, $(a,b)\neq (0,0)$), where $f(z)$ is a
meromorphic function with quasi-periodicity,
\begin{equation}
    f(z+1) = \eps^a f(z), \qquad f(z+\tau) = \eps^b f(z).
\label{period:gtw:ell}
\end{equation}

For each point $P$ on $E$, we define a Lie algebra,
\begin{equation}
     \g^P := (\gtw \tensor_{\calO_E} \K_E)^{\wedge}_P,
 \label{def:loop}
\end{equation}
where $\K_E$ is the sheaf of meromorphic functions on $E$ and
$(\cdot)_P^{\wedge}$ means the completion of the stalk at $P$ with
respect to the maximal ideal $\m_P$ of $\calO_{E,P}$, the stalk of the
structure sheaf. The Lie algebra $\g^P$ is (non-canonically) isomorphic
to the loop Lie algebra $\g((z-z_0))$, where $z_0$ is the coordinate of
$P$. The subspace
\begin{equation}
    \g^P_+:=(\gtw)_P^{\wedge}\simeqq\g[[z-z_0]]
\label{def:loop+}
\end{equation}
of $\g^P$ is a Lie subalgebra.

Let us fix mutually distinct points $Q_1, \ldots, Q_L$ on $E$ whose
coordinates are $z = z_1, \ldots, z_L$ and put $D:=\{Q_1,\dots,Q_L\}$.
We shall also regard $D$ as a divisor on $E$ (i.e.,
$D=Q_1+\cdots+Q_L$). The Lie algebra $\g^D:=\bigoplus_{i=1}^L\g^{Q_i}$
has a 2-cocycle defined by
\begin{equation}
    \caff(A,B) := \sum_{i=1}^L \caff_{,i}(A_i,B_i),\qquad
    \caff_{,i}(A_i,B_i) := \Res_{Q_i} (d A_i | B_i),
\label{def:cocycle}
\end{equation}
where $A=(A_i)_{i=1}^L, B=(B_i)_{i=1}^L \in \g^D$, $\Res_{Q_i}$ is the
residue at $Q_i$ and $d$ is the exterior derivation. (The symbol
``$\caff$'' stands for ``Cocycle defining the Affine Lie algebra''.)  We
denote the central extension of $\g^D$ with respect to this cocycle by
$\hat\g^D$:
\begin{equation}
    \hat\g^D := \g^D \oplus \C \khat,
\label{affine}
\end{equation}
where $\khat$ is a central element. Explicitly the bracket of $\hat\g^D$
is represented as
\begin{equation}
  [A, B] = ([A_i, B_i]^\circ)_{i=1}^L \oplus \caff(A,B) \khat
  \quad
  \text{for $A,B\in\g^D$,}
\label{def:aff-alg-str}  
\end{equation}
where $[A_i, B_i]^\circ$ are the natural bracket in $\g^{Q_i}$.  The Lie
algebra $\hat\g^P$ for a point $P$ is nothing but the affine Lie algebra
$\hat\g$ of type $A^{(1)}_{N-1}$ (a central extension of the loop
algebra $\g((t-z))=sl_N\bigl(\C((t-z))\bigr)$).  

The affine Lie algebra $\hat\g^{Q_i}$ can be regarded as a subalgebra of
$\hat\g^D$.  The subalgebra $\g^{Q_i}_+$ of $\g^{Q_i}$ (cf.\
\eqref{def:loop+}) can be also regarded as a subalgebra of
$\hat\g^{Q_i}$ and $\hat\g^D$.

Let $\gout$ be the space of global meromorphic sections of $\gtw$ which are
holomorphic on $E$ except at $D$:
\begin{equation}
    \gout := \Gamma(E, \gtw(\ast D)).
\label{def:gout:ell}
\end{equation}
The residue theorem implies that we can regard $\gout$ as a Lie
subalgebra of $\hat\g^D$ by mapping an element of $\gout$ to its germs
at $Q_i$'s.

\begin{definition}
\label{def:CC,CB:ell}
The space of {\em conformal coinvariants} $\CC_E(M)$ and that of {\em
conformal blocks} $\CB_E(M)$ over $E$ associated to
$\hat\g^{Q_i}$-modules $M_i$ with the same level $\khat = k$ are defined
by
\begin{equation}
    \CC_E(M) := M/\gout M, \quad
    \CB_E(M) := \Hom_\C(M/\gout M, \C),
\label{def:cc,cb:ell}
\end{equation}
 where $M := \bigotimes_{i=1}^L M_i$. The module $M_i$ is referred to as
 ``a module inserted at the point $Q_i$''.

\end{definition}

\section{Family of elliptic curves}
\label{sec:family}

In this section we construct a family of elliptic curves $\calE$ with a
singular fibre and its covering $\tilde\calE$. A twisted Lie algebra
bundle $\gtw$ over the generic fibres of $\calE$ is one of the main
representation theoretical data in the twisted WZW model but it does not
directly extends to the singular fibre. Hence we pull it back to
$\tilde\calE$ and trivialise it. The sections of $\gtw$ on the singular
fibre is understood as sections of the trivial bundle on $\tilde\calE$
invariant under the action of the covering transformation group.

The construction of $\tilde\calE$ is almost the same as that of the
analytic fibre space of elliptic curves in \cite{wol:83}. We use $N$
patches $U_k$ ($k\in \Z/N\Z$):
\begin{equation}
    U_k := 
   \{(q,x_k,y_k) \mid 
     |q|<1, x_k y_k = q, |x_k| < |q|^{-1}, |y_k| < |q|^{-1} \}.
\label{def:Uk}
\end{equation}
We denote $(q,x,y) \in U_k$ by $(q,x,y)_k$. The universal curve
$\tilde\calE$ is defined by
\begin{equation}
    \tilde\calE := \Bigl. \bigsqcup_{k\in\Z/N\Z} U_k \Bigr/ \sim,
\label{def:tildeC}
\end{equation}
where the equivalence relation $\sim$ is defined by
\begin{equation}
    (q,x_k,y_k)_k \sim (q,x'_{k+1},y'_{k+1})_{k+1} \text{\ when\ }
    x_k y'_{k+1} = 1.
\label{gluing}
\end{equation}
We have an analytic fibre space $\tilde\pi: \tilde\calE \to \Delta$ over
$\Delta = \{ q \mid |q|<1\}$. The fibre over $q \neq 0$ is an elliptic
curve $\C^{\times}/q^{N\Z}$ and the fibre over $q=0$ is singular with
ordinary double points $(0,0,0)_k \in U_k$ ($k\in \Z/N\Z$).

Let $C_N = \{\bar{0},\bar{1},\dots,\overline{N-1}\}$ be the cyclic group of
order $N$. The group $C_N^2 = C_N \times C_N$ acts on $\tilde\calE$ from
the right as follows: it is enough to define the action of the
generators $(1,0)$ and $(0,1)$ of $C_N^2$.
\begin{equation}
    (q,x_k,y_k)_k \cdot (1,0) = (q,\eps^{-1} x_k, \eps y_k)_k, \qquad
    (q,x_k,y_k)_k \cdot (0,1) = (q,x_k, y_k)_{k-1}.
\label{CN-action}
\end{equation}
The universal curve $\calE$ is defined set-theoretically as the quotient
space of $\tilde\calE$ by this action:
\begin{equation}
    \calE := \tilde\calE/C_N^2.
\label{def:C}
\end{equation}
The canonical projection to $\Delta$ is denoted by $\pi: \calE \to
\Delta$. The fibre $\pi^{-1}(0)$ is the singular curve with one ordinary
double point.

\begin{figure}[h]
 \begin{center}
  \psfrag{tildeC}{$\tilde\calE$}
  \psfrag{C}{$\calE$}
  \psfrag{q=/0}{$q\neq0$}
  \psfrag{q=0}{$q=0$}
  \includegraphics{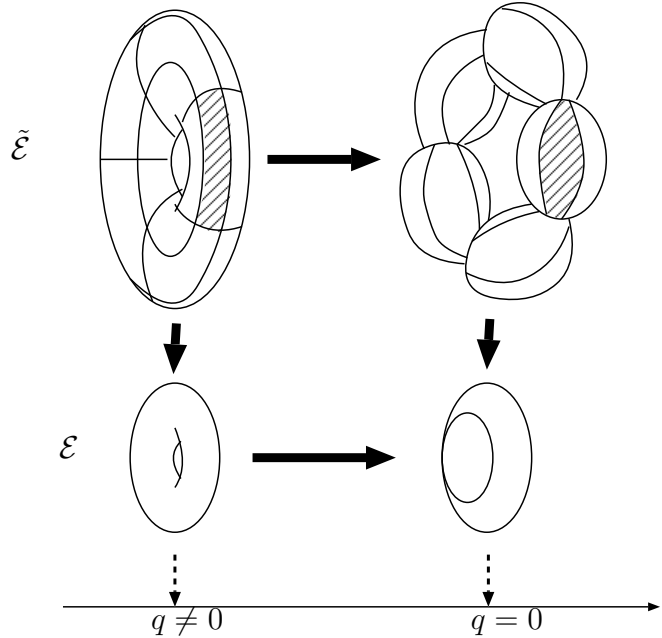}
 \end{center}
\caption{Degeneration of an elliptic curve and its $N^2$-covering.}
\end{figure}

The family of $L$-pointed elliptic curves and its covering are defined
by
\begin{equation}
    \X := \calE \times_{\Delta} S, \qquad
    \tilde\X := \tilde\calE \times_{\Delta} S.
\label{def:X,tildeX}
\end{equation}
Here the base space $S$ is the fibre product of $\calE$'s
without diagonals:
\begin{equation}
 \begin{split}
    S := &\{(q; Q_1,\dots,Q_L) \in 
    \overbrace{\calE \times_\Delta \cdots \times_\Delta \calE}^L
\\
    &\mid
    Q_i \neq Q_j \, (i\neq j), Q_i \neq [(0,0,0)_k] \text{\ for any k.}
    \},
 \end{split}
\label{def:S}
\end{equation}
where $Q_i$ is the point of the $i$-th $\calE$ in the fibre product with
$\pi(Q_i) = q$. We exclude the degeneration of the types $Q_i \to Q_j$
and $Q_i \to$ (node). We denote the canonical projections $\X\to S$,
$\tilde\X \to S$ and $\tilde\X \to \X$ by $\pi_{\X/S}$,
$\pi_{\tilde\X/S}$ and $\pi_{\tilde\X/\X}$ respectively. The fibres of
$\tilde\X$ and $\X$ over $S_0 := \{(0;Q_1,\dots,Q_L) \in S\}$ are
singular curves.

The section $q_i$ of $\X \to S$ is defined by
\begin{equation}
    q_i ((q; Q_1,\dots,Q_L)) := (q; Q_i; Q_1, \dots, Q_L).
\label{def:si}
\end{equation}
We denote the divisor $[q_1(S)] + \cdots + [q_L(S)]$ by $D$.

The group $C_N^2$ acts on $\tilde\X$ naturally as covering
transformation and on $\g$ by
\begin{equation}
    (m,n) \cdot A = (\gamma^m \beta^n) A (\gamma^m \beta^n)^{-1}.
\label{CN2-on-g}
\end{equation}
The twisted Lie algebra bundle $\gtwX$ on $\overset{\circ}{\X} := \X
\setminus \{\text{singular points}\}$ is defined as the associated
bundle to the $C_N^2$-principal bundle $\overset{\circ}{\tilde\X}:=
\tilde\X \setminus\{\text{singular points}\} \to \overset{\circ}{\X}$:
\begin{equation}
    \gtwX := \overset{\circ}{\tilde\X} \times_{C_N^2} \g.
\label{def:gtwX}
\end{equation}
It is obvious that the restriction of $\gtwX$ to a fibre of $\X$ at a
point $(q; Q_1,\dots,Q_L)$ ($q\neq 0$) is the bundle $\gtw$ on the
elliptic curve $\C^{\times}/q^\Z$ defined by \eqref{def:gtw}.

Sheaf version of affine Lie algebras $\g^P$, $\g^P_+$, $\g^D$ and
$\hat\g^D$ (cf.\ \eqref{def:loop}, \eqref{def:loop+}, \eqref{affine})
are $\calO_S$-Lie algebras defined by
\begin{equation}
\begin{aligned}
  & \g_S^{Q_i} 
  := \pi_{\X/S,\ast} (\gtwX(\ast Q_i))_{Q_i}^{\wedge},
  \quad &
  & \g_{S,+}^{Q_i}
  := \pi_{\X/S,\ast} (\gtwX)_{Q_i}^{\wedge},
\\
  & \g_S^D
  := \pi_{\X/S,\ast} (\gtwX(\ast D))_D^{\wedge}
  = \bigoplus_{i=1}^L \g_S^{Q_i},
  \quad &
  & \g_{S,+}^D
  := \pi_{\X/S,\ast} (\gtwX)_D^{\wedge}
  = \bigoplus_{i=1}^L \g_{S,+}^{Q_i}.
\end{aligned}
\label{def:Os-Lie-alg}
\end{equation}
Since we assume that sections $Q_i$ do not touch the singular point
of the singular fibre, the definitions are the same as those for the
non-singular case, (3.11) of \cite{kur-tak:97}. The central extensions
of $\g^D_S$ and $\g^D_{S,+}$ are defined by the cocycle
\eqref{def:cocycle} with the coefficients in $\calO_S$:
\begin{equation}
    \hat\g^D_S := \g^D_S \oplus \calO_S \hat k, \qquad
    \hat\g^D_{S,+} := \g^D_{S,+} \oplus \calO_S \hat k.
\label{def:affine/S}
\end{equation}

The Lie subalgebra of meromorphic sections $\gout \subset \hat\g^D$,
\eqref{def:gout:ell}, would be replaced by $\pi_{\X/S,\ast}(\gtwX(\ast
D))$ if there were no singularity, as was the case in
\cite{kur-tak:97}. Taking the singular fibre into account, we modify
this naive definition as follows:
\begin{equation}
    \goutX = 
    \bigl(\pi_{\tilde\X/S,\ast}(\g\tensor\calO_{\tilde\X} (\ast\tilde D))
    \bigr)^{C_N^2},
\label{def:goutX}
\end{equation}
where $\tilde D$ is the $C_N^2$-orbit of the divisor $D$ and
$(\cdot)^{C_N^2}$ denotes the $C_N^2$-invariant section of the
equivariant locally free sheaf $\g\tensor \calO_{\tilde\X}$. In other
words, a $\g$-valued meromorphic function $f(s,P)$ ($s =
(q;Q_1,\dots,Q_L) \in S$, $P \in \tilde\calE|_s$) belongs to $\goutX$ if
and only if it satisfies
\begin{equation}
    f(s,(1,0)\cdot P) = \Ad\gamma(f(s,P)), \qquad
    f(s,(0,1)\cdot P) = \Ad\beta(f(s,P)),
\label{sec-goutX}
\end{equation}
where $(m,n) \cdot$ ($m,n \in\Z$) is the left action of the generators
of $C_N^2$ on the fibre of $\tilde\calE$ defined by $(m,n)\cdot P := P
\cdot (-m,-n)$. (See \eqref{CN-action}.)

Our construction is so explicit that we have an explicit basis of
$\goutX$. In \cite{kur-tak:97} meromorphic functions $w_{ab}(\tau;t)$
($a,b = 0,\dots,N-1$, $\Im\tau>0$, $t\in\C$) characterised by the
following properties were introduced:
\begin{itemize}
 \item Additive quasi-periodicity: $w_{ab}(\tau;t+1) = \eps^a
       w_{ab}(\tau;t)$, $w_{ab}(\tau;t+\tau) = \eps^b w_{ab}(\tau;t)$;
 \item As a function of $t\in\C$, $w_{ab}(\tau;t)$ has a simple pole
       with residue $1$ at $\Z + \Z\tau$.
\end{itemize}
Let us denote this function by $\wadd_{ab}(\tau;t)$. (The superscript
``add'' stands for ``additive''.) Let us rewrite it to a
multiplicatively quasi-periodic function $\wmul_{ab}(q;u)$ as follows:
\begin{multline}
    \wmul_{ab}(q;u)
    :=
    \frac{2\pi i u^a}{u^N - 1} \times
\\
    \times
    \frac{ (q^{N-a} \eps^b    u^{-N}; q^N)_\infty 
           (q^a     \eps^{-b} u^N   ; q^N)_\infty}
         { (q^{N-a} \eps^b          ; q^N)_\infty
           (q^N               u^{-n}; q^N)_\infty
           (q^a     \eps^{-b}       ; q^N)_\infty
           (q^N               u^N   ; q^N)_\infty},
\label{def:wmul}
\end{multline}
where $(x;q)_\infty = \prod_{n=0}^\infty (1-x q^n)$ is the standard
infinite product symbol. The function $\wmul_{ab}(q;u)$ is related to
$\wadd_{ab}$ by $\wmul_{ab}(q;u) = \wadd_{ab}(N \log u/2\pi i, N \log
q/2\pi i)$ when $q\neq 0$, that is, we replaced the arguments of
$\wadd_{ab}$ by $e^{2\pi i z} = u^N$, $e^{2\pi i \tau} = q^N$ and used
the product formula for the theta function. When $q=0$, $\wmul_{ab}$
becomes a rational function of $u$:
\begin{equation}
    \wmul_{ab}(0;u) =
    \begin{cases}
      2\pi i u^a (u^N - 1)^{-1}, & (a\neq 0),\\
      2\pi i(1-\eps^b)^{-1} (u^N - \eps^b)(u^N - 1)^{-1}, &(a=0).
    \end{cases}
\label{wmul(q=0)}
\end{equation}
The important property of $\wmul_{ab}$ is that it inherits the
quasi-periodicity of $\wadd_{ab}$:
\begin{equation}
    \wmul_{ab}(q;\eps u) = \eps^a \wmul_{ab}(q;u), \qquad
    \wmul_{ab}(q;   q u) = \eps^b \wmul_{ab}(q;u).
\label{wmul:periodicity}
\end{equation}
We define the function $w_{ab,i}(P)$ ($i=1,\dots,L$) on $\tilde\X$ in
terms of $\wmul_{ab}$: for $(q,x_k,y_k)_k \in U_k$ (cf.\
\eqref{def:Uk}), $q\neq 0$, 
\begin{equation}
    w_{ab,i}(((q,x_k,y_k)_k, s)) := \wmul_{ab}(q; q^{k-k'}x_k/x'_{k'}).
\label{def:wab}
\end{equation}
Here $s = (q; Q_1,\dots,Q_L) \in S$ and we fix an index $k'$ to express
$Q_i$ as a point $(q,x'_{k'},y'_{k'})_{k'}$ in $U_{k'}$. (The function
$w_{ab,i}$ is determined up to this choice.) This function is extended
to the points with $q=0$. (For example, $w_{ab,i}((0,x_0,y_0)_0)$ is the
rational function \eqref{wmul(q=0)} of $x=x_0/x'_{0}$ when $Q_i$ is
represented as a point $(0; x'_0, y'_0)_0$ in $U_0$.) The main
properties of this function are
\begin{itemize}
 \item Quasi-periodicity with respect to $C_N^2$-action:
 \begin{equation}
     w_{ab,i}((1,0)\cdot P) = \eps^a w_{ab,i}(P), 
\qquad
     w_{ab,i}((0,1)\cdot P) = \eps^b w_{ab,i}(P),
\label{wab:periodicity}
 \end{equation}

 \item All poles are simple and located at $Q_i$ modulo $C_N^2$
       action. 
\end{itemize}

It is easy to see that any section of $\goutX$ is a linear combination
of $J_{ab} \tensor w_{ab,i}(P)$'s and their derivatives along the fibre. 

\begin{lemma}
\label{lem:gD=gout+regular}
 $\hat\g_S^D = \goutX \oplus \hat\g_{S,+}^D$.
\end{lemma}

\begin{proof}
 The singular part of an element of $\hat\g_S^D$ can be expressed by a
 linear combination of $J_{ab} \tensor w_{ab,i}(P)$ and derivatives in a
 unique way. Subtracting such combination which belongs to $\goutX$,
 we end up with a regular element of $\hat\g_{S,+}^D$.
\end{proof}

\section{Sheaves of conformal coinvariants and conformal blocks}
\label{sec:cc,cb}

In this section we introduce the sheaf $\sCC$ of conformal coinvariants
and the sheaf $\sCB$ of conformal blocks and show their basic
properties.

Definitions of $\sCC$ and $\sCB$ are literally the same as those for the
non-singular case, Definition 3.3 of \cite{kur-tak:97}.
\begin{definition}
\label{def:sCC-sCB} For any $\hat\g_S^D$-module $\M$ of level $k$ (i.e.,
 $\khat$ acts as $k\cdot\id$), we define the {\em sheaf $\sCC(\M)$ of
 conformal coinvariants} and the {\em sheaf $\sCB(\M)$ of conformal
 blocks} by
\begin{align}
  \sCC(\M)
  & := \M/\goutX\M,
  \label{def:sheaf-conf-coinvariant}
\\
  \sCB(\M)
  & := \HOM_{\calO_S}(\sCC(\M), \calO_S)
  \label{def:sheaf-conf-block}.
\end{align}
We can regard $\sCC(\cdot)$ as a covariant right exact functor from the
category of $\hat\g_S^D$-modules to that of $\calO_S$-modules and similarly
$\sCB(\cdot)$ as a contravariant left exact functor.
\end{definition}

The goal of this paper is to prove that $\sCC$ and $\sCB$ are locally
free. Since $\sCB$ is the dual of $\sCC$, we mainly discuss about $\sCC$
and briefly mention on $\sCB$ when it is necessary.

We assume that the $\hat\g_S^D$-module $\M$ are of the following type:
\begin{equation}
    \M = \calO_S \tensor \bigotimes_{i=1}^L M_i
       = \bigotimes_{i=1}^L (\calO_S \tensor M_i),
\label{M:finiteness}
\end{equation}
where each $M_i$ is a quotient of a $\hat\g$-Weyl module $M(V_i):=
\Ind_{\hat\g^{Q_i}_+}^{\g^{Q_i}\oplus\C\hat k} V_i$ of level $k$ for a
finite dimensional irreducible $\g$-module $V_i$. (See \S6 of
\cite{tak:04} or \S2.4 of \cite{kaz-lus:93}.) To endow $\M$ in
\eqref{M:finiteness} with the $\hat\g_S^D$-module structure, we need to
fix the coordinate of $\calE$ and the trivialisation of $\gtwX$, which
is irrelevant to the statements of theorems below. In the concrete
computations, we use the coordinates and the trivialisation obtained
naturally from the construction in \secref{sec:family}.

\begin{proposition}
\label{prop:coherence}
 $\sCC(\M)$ is a coherent $\calO_S$-sheaf. When all $M_i$'s are Weyl
 modules, $\sCC(\M)$ is (and hence $\sCB(\M)$ is) locally free.
\end{proposition}

\begin{proof}
 \lemref{lem:gD=gout+regular} makes it possible to apply the same
 argument as the proof for the non-singular case, Corollary 3.5 of
 \cite{kur-tak:97}. In fact, if each $M_i$ in \eqref{M:finiteness} is
 a Weyl module $M(V_i)$, $\M$ is expressed as
\begin{equation}
    \M 
    = U_S(\hat\g^D_S) \tensor_{U_S(\hat\g^D_{S,+})} (V\otimes\calO_S)
    \isoot U_S(\goutX) \tensor_{\calO_S} (V\otimes\calO_S),
\label{M=U(goutX)*V}
\end{equation}
 by the Poincar\'e-Birkhoff-Witt theorem ($V = \tensor_{i=1}^L
 V_i$). Here $U_S(\cdot)$ denotes the universal $\calO_S$-enveloping
 algebra.  Modding out by $\goutX \M$, we have
\begin{equation}
    \sCC(\M) = \M/\goutX\M \isoot V \tensor \calO_S,
\label{sCC(Weyl)=free}
\end{equation}
 which means that $\sCC(\M)$ is locally free, in particular,
 coherent.

 If $M_i$'s are quotients of Weyl modules, $\M$ is a quotient of
 $\bigotimes M(V_i) \tensor \calO_S$. Hence by the right exactness of the
 functor $\sCC$, $\sCC(\M)$ is a quotient of a coherent sheaf, and
 therefore coherent.
\end{proof}

In \secref{sec:loc-free} we prove the locally freeness of $\sCC(\M)$ for
integrable $M_i$'s, examining the behaviour of $\sCC(\M)$ at the
boundary of the moduli space ($S_0 = \{q=0\} \subset S$) carefully.

\section{Sheaf version of trigonometric and orbifold
 WZW model } 
\label{sec:sheaf-tri/orb}

Everything in previous two sections can be restricted on $S_0$, namely
on the configuration space of points on a singular rational curve with
one ordinary double point. (As we have mentioned, the restriction of the
functions $w_{ab,i}$ needs special care.) Hence we can define the
corresponding sheaves $\sCC(\M)$ and $\sCB(\M)$ which we denote by
$\sCCtri(\M)$ and $\sCBtri(\M)$. The subscript ``trig'' is put here
because, as we shall see below, there are connections on them expressed
in terms of the trigonometric $r$ matrix.

In the proof of \thmref{thm:loc-free} we shall use the sheaf of
conformal coinvariants of the orbifold WZW model,
$\sCCorb(\M_0\tensor\M\tensor\M_\infty)$ on $S_0$, where $\M_* = M_*
\tensor \calO_{S_0}$ for a $\hat\g^{(\ast)}$-module $M_\ast$ ($\ast =
0,\infty$). We shall recall the definition of the twisted affine
algebras $\hat\g^{(0)}$ and $\hat\g^{(\infty)}$ and the details of
$\CC_\orb$ soon later. Here we only say that $\sCCorb$ is defined
exactly in the same way as $\sCCtri(\M)$ if we replace the degenerate
elliptic curve (the fibre of $\calE$ at $q=0$) by the orbifold
$\Proj^1/C_N$. Note that $S_0$ can be regarded as the configuration
space of points on $\Proj^1/C_N$. (See \S3 of \cite{tak:04}.)

Let $\M$ be as in \eqref{M:finiteness} and $M_*$ be a
$\hat\g^{(*)}$-Verma module. Then \propref{prop:coherence} holds for
$\sCCtri$, $\sCBtri$, $\sCCorb$ and $\sCBorb$ as well. In fact, we can
prove locally freeness under this assumption.

\begin{proposition}
\label{prop:loc-free:trig/orb}
 (i)
 The sheaf $\sCCtri(\M)$ and the sheaf $\sCBtri(\M)$ are locally
 free $\calO_{S_0}$-sheaves.

 (ii)
 The sheaf $\sCCorb(\M_0\tensor\M\tensor\M_\infty)$ and the sheaf
 $\sCBorb(\M_0\tensor\M\tensor\M_\infty)$ are locally free
 $\calO_{S_0}$-sheaves.
\end{proposition}

\begin{proof}
 (i) 
 The proof of the locally freeness for the non-singular case, Corollary
 5.3 of \cite{kur-tak:97} is true also in this case: $\sCCtri(\M)$ is
 coherent as shown at the end of \secref{sec:cc,cb} and there is a
 connection and $D$-module structure on it, which implies that it is
 locally free $\calO_{S_0}$-sheaf. The only difference is that we do not
 change the curve itself (the modulus $q$ is fixed to $0$) in the
 present case, and hence there is nothing corresponding to the
 connection in the direction of $\partial/\partial\tau$ in
 \cite{kur-tak:97}. The connection in the direction of
 $\partial/\partial z_i$ ($z_i$ is the coordinate of $Q_i$) is
 $\nabla_i = \partial/\partial z_i - \rho_i(T[-1])$ (cf.\ (5.14) of
 \cite{kur-tak:97}) as is well-known, where $\rho_i$ is the
 representation of the Virasoro algebra on $M_i$ constructed via the
 Sugawara construction and $T[-1]$ is one of the Virasoro generator,
 usually denoted by $L_{-1}$.

 (ii)
 We might proceed as the proof of (i) from the beginning but the short
 cut is to use the result of (i). The coherence of $\sCCorb(\M_0 \tensor
 \M \tensor \M_\infty)$ having been proved, we have only to check that
 the above connection operators $\nabla_i$ ($i=1,\dots,L$) on
 $\sCCtri(\M)$ also define the flat connection on $\sCCorb$. What we
 need to check is
\begin{itemize}
 \item $[\nabla_i, \gout^\orb] \subset \gout^\orb$,
 \item $[\nabla_i, \nabla_j] = 0$ on $\sCCorb$,
\end{itemize}
 which is proved in the same way as in the case of the ordinary WZW
 model, e.g., Lemma 4 of \cite{f-f-r:94}.
\end{proof}

The connection on $\sCCtri(\M)$ mentioned in the proof of (i) is
obtained by the degeneration $q \to 0$ of the elliptic
Knizhnik-Zamolodchikov connection in \cite{eti:94} and
\cite{kur-tak:97}. They are expressed as the first order differential
operators on $V \tensor \calO_{S_0}$ in terms of the {\em trigonometric $r$
matrix}. In fact, by tracing the argument which leads to the explicit
form (Theorem 5.9 in \cite{kur-tak:97}) of the connection, we have only
to replace the functions $w_{ab}(z_j-z_i)$ there with $w_{ab,i}(Q_j)$
which is expressed by the rational function of the form
\eqref{wmul(q=0)} on $S_0$. Hence the KZ equation for the WZW model on
the degenerate elliptic curve is the {\em trigonometric KZ equation\/}.

\begin{lemma}
\label{lem:CCtri/orb-fibre}
 (i)
 The fibre of $\sCCtri(\M)$ at $s\in S_0$, $\sCCtri(\M)|_s$, is
 isomorphic to $\CC_\trig(M)$, the space of conformal coinvariants of
 the trigonometric WZW model for the geometric data corresponding to
 $s$.

 (ii)
 The fibre of $\sCCorb(\M)$ at $s$, $\sCCorb(\M)|_s$, is isomorphic to
 $\CC_\orb(M)$, the space of conformal coinvariants of the
 orbifold WZW model for the geometric data corresponding to $s$. 
\end{lemma}

See Definition 3.2 of \cite{tak:04} for the definition of $\CC_\trig$
and $\CC_\orb$.

\begin{proof}
 We can modify the proof of the corresponding statement for the
 non-singular case, Corollary 5.4 of \cite{kur-tak:97}. For example, for
 the case (i), the isomorphism $\goutX|_{s} \simeqq \gout^\trig$ is a
 consequence of the existence of the sections $w_{ab,i}$ and their
 derivatives (cf.\ the end of \secref{sec:family}) which span both
 $\goutX|_s$ and $\gout^\trig$. The rest of the proof can be translated
 to the present case without change.

 The proof of (ii) is similar.
\end{proof}

Combining \propref{prop:loc-free:trig/orb}, \lemref{lem:CCtri/orb-fibre}
and Theorem 5.1 of \cite{tak:04}, we have the following isomorphism:
\begin{equation}
    \iota: \sCCtri(\M) \isoto
    \bigoplus_{\lambda \in \wt(V)}
    \sCCorb(\M^{(0)}_{\tilde\lambda} \tensor \M 
    \tensor \M^{(\infty)}_{\tilde\lambda'}),
\label{factor:trig->orb}
\end{equation}
where $V = \bigotimes_{i=1}^L V_i$ is the $\g$-module generating $M =
\bigotimes_{i=1}^L M_i$ (recall $M_i$ is a quotient of Weyl module
$M(V_i)$), $\wt(V)$ is the set of its weights, 
\begin{equation}
 \begin{split}
    \tilde \lambda  &:= - \lambda \circ (1- \Ad\beta^{-1})^{-1}
    = \lambda \circ (1- \Ad\beta)^{-1} \circ \Ad\beta,
\\
    \tilde \lambda' &:= - \lambda \circ (1- \Ad\beta)^{-1},
 \end{split}
\label{def:tilde-lambda}
\end{equation}
$\M^{(*)}_\mu = M^{(*)}_\mu \tensor \calO_{S_0}$ for a
Verma module $M^{(*)}_\mu$ of $\hat\g^{(*)}$ with the highest weight
$\mu$ (cf.\ Definition 4.1 (i) of \cite{tak:04}).

In \secref{sec:loc-free} the modules $M_i$ are assumed to be integrable
highest weight modules. (cf.\ Chapter 10 of \cite{kac:90}.)  In this
case the above result can be refined. For this purpose we recall the
details of the orbifold WZW model defined in \S3 of \cite{tak:04}.

Let us denote the standard coordinate of $\Proj^1(\C)$ by $t$. The
cyclic group $C_N$ acts as $t \mapsto \eps^a t$ ($a \in C_N$) and the
quotient $E_\orb = \Proj^1/C_N$ is an ordinary orbifold.

The definition of the space of conformal coinvariants/blocks of the
orbifold WZW model on $E_\orb$ is almost the same as that on elliptic
curves, \defref{def:CC,CB:ell}, except that we also insert modules to
the singular points $0$ and $\infty$. The Lie algebra $\gout$ in
\eqref{def:cc,cb:ell} is replaced by $\gout^\orb$ which consists of
$\g$-valued meromorphic functions $f(t)$ on $\Proj^1$ such that:
(1) poles belong to $\{0, Q_1,\dots, Q_L, \infty\}$;
(2) $f(\eps t) = \Ad\gamma (f(t))$.
Accordingly, the module inserted at $0$ is the $\hat\g^{(0)}$-module and
the module inserted at $\infty$ is the $\hat\g^{(\infty)}$-module, where
\begin{align}
    \hat\g^{(0)} &= 
    \bigoplus_{\substack{a,b=0,\dots,N-1\\ (a,b)\ne(0,0)}}
    \bigoplus_{m\in\Z} 
    \C J_{a,b} \tensor t^{a+mN}
    \oplus \C\khat,
\label{g(0)}
\\
    \g^{(\infty)} &= 
    \bigoplus_{\substack{a,b=0,\dots,N-1\\ (a,b)\ne(0,0)}}
    \bigoplus_{m\in\Z} 
    \C J_{a,b} \tensor t^{a+mN}
    \oplus \C\khat
\label{g(infty)}
\\
    &= \bigoplus_{\substack{a,b=0,\dots,N-1\\ (a,b)\ne(0,0)}}
    \bigoplus_{m\in\Z} 
    \C J_{a,b} \tensor s^{-a+mN}
    \oplus \C\khat.
    \qquad (s = t^{-1})
\nonumber
\end{align}
The cocycles which defines the central extension of $\g^{(0)}$ and
$\g^{(\infty)}$ are:
\begin{equation}
    \caff_{,0}(A,B):= \frac{1}{N}\Res_{t=0} (d A|B), \qquad
    \caff_{,\infty}(A,B):= \frac{1}{N}\Res_{s=0} (d A|B),
\label{def:caff:0,infty}
\end{equation}
for $A,B\in\g^{(0)}$ and $A,B\in\g^{(\infty)}$ respectively. ($\g^{(*)}$
is the loop algebra part of $\hat\g^{(*)}$.)

As Etingof showed (Lemma 1.1 of \cite{eti:94}), $\hat\g^{(0)}$ and
$\hat\g^{(\infty)}$ are isomorphic to the ordinary affine Lie algebra
$\hat\g$ of $A^{(1)}_{N-1}$ type. Explicitly the isomorphism $\phi_0:
\hat\g^{(0)} \isoto \hat\g$ is defined by
\begin{equation}
 \begin{gathered}
    \phi_0(E_{i,i+1}\tensor t) = e_i, \qquad
    \phi_0(E_{i+1,i}\tensor t^{-1}) = f_i,
\\
    \phi_0(E_{N,1} \tensor t) = e_0, \qquad
    \phi_0(E_{1,N} \tensor t^{-1}) = f_0,
\\
    \phi_0(H_{i,i+1} \tensor 1) = \alpha_i^{\vee} - \frac{\khat}{N}, \qquad
    \phi_0(\khat) = \khat,
 \end{gathered}
\label{def:phi0}
\end{equation}
for $i=1,\dots,N-1$, where $E_{ij}$ is the matrix unit, $H_{ij} = E_{ii}
- E_{jj}$, $e_i$, $f_i$ ($i=0,\dots,N-1$) are the Chevalley generators
of $\hat\g$ and $\alpha_i^{\vee}$ ($i=0,\dots,N-1$) are coroots of
$\hat\g$. (cf.\ \S6.2 and \S7.4 of \cite{kac:90}.) Since the positive
powers of $s$ (cf.\ \eqref{g(infty)}) kill the highest weight vector of
the Verma module $M^{(\infty)}_{\mu}$ inserted at $\infty$ (cf.\
Definition 4.1 (i) of \cite{tak:04}), we identify $\hat\g^{(\infty)}$
and $\hat\g$ through an isomorphism which is essentially a composition
of $\phi_0$ and the Chevalley involution:
\begin{equation}
 \begin{gathered}
    \phi_{\infty}(E_{i,i+1}\tensor s^{-1}) = - f_i, \qquad
    \phi_{\infty}(E_{i+1,i}\tensor s) = - e_i,
\\
    \phi_{\infty}(E_{N,1} \tensor s^{-1}) = - f_0, \qquad
    \phi_{\infty}(E_{1,N} \tensor s) = - e_0,
\\
    \phi_{\infty}(H_{i,i+1} \tensor 1) 
    = - \alpha_i^{\vee} + \frac{\khat}{N}, \qquad
    \phi_{\infty}(\khat) = \khat,
 \end{gathered}
\label{def:phi00}
\end{equation}
for $i=1,\dots,N-1$. Identified through $\phi_0$ and $\phi_\infty$, the
Verma modules of $\hat\g^{(0)}$ and $\hat\g^{(\infty)}$ are Verma
modules of $\hat\g$ in ordinary sense. 

In this section, $e_i$ and $f_i$ denote the Chevalley generators of
$\hat\g$ identified with the elements in $\hat\g^{(0)}$ or
$\hat\g^{(\infty)}$ by means of $\phi_0$ or $\phi_\infty$.

\begin{proposition}
\label{prop:ccorb:integrable} 
 Assume that all $M_i$ ($i=1,\dots,L$) are integrable highest weight
 modules and that $M^{(\ast)}_{\mu_\ast}$ ($\ast = 0, \infty$) is a
 Verma modules of $\hat\g^{(\ast)}$.

 Then $\CC_\orb(M^{(0)}_{\mu_0} \tensor M \tensor
 M^{(\infty)}_{\mu_\infty})$ is $0$ unless $\mu_0$ and $\mu_\infty$ are
 dominant integral weights of $\hat\g$ identified with $\hat\g^{(0)}$
 and $\hat\g^{(\infty)}$. If $\mu_0$ and $\mu_\infty$ are dominant
 integral weights, 
\begin{equation}
    \CC_\orb(M^{(0)}_{\mu_0}\tensor M\tensor M^{(\infty)}_{\mu_\infty}) 
    \simeqq
    \CC_\orb(L^{(0)}_{\mu_0}\tensor M\tensor L^{(\infty)}_{\mu_\infty}),
\label{ccorb:integrable}
\end{equation} 
 where $L^{(\ast)}_{\mu_*}$ ($\ast = 0, \infty$) is the irreducible
 quotient of $M^{(\ast)}_{\mu_\ast}$.
\end{proposition}

\begin{rem}
 In physics context, this proposition is a consequence of the
 propagation of the null field. See \S4 of \cite{zam:89}. The author
 thanks Yasuhiko Yamada for this comment.
\end{rem}

\begin{rem}
 \propref{prop:ccorb:integrable} is in sharp contrast to the Weyl module
 case. See \S6 of \cite{tak:04}.
\end{rem}

\begin{proof}
 The following lemma shall be proved later.

\begin{lemma}
\label{lem:cc-reduction}
 Let $N_\ast$ ($\ast = 0,\infty$) be a quotient of the Verma module
 $M^{(\ast)}_{\mu_\ast}$. Suppose $v_\kappa\in N_0$ is a singular vector
 of weight $\kappa$ which is not a dominant integral weight. (The weight
 $\kappa$ may possibly be the highest weight $\mu_0$.) Then for any $v
 \in M$ and $v_\infty \in N_\infty$,
\begin{equation}
    v_\kappa \tensor v \tensor v_\infty \equiv 0
    \mod{\gout^\orb (N_0 \tensor M \tensor N_\infty)}.
\label{cc-reduction}
\end{equation}
 The same is true for a singular vector $v_\kappa \in N_\infty$.
\end{lemma}

 The first statement of \propref{prop:ccorb:integrable} is a consequence
 of \lemref{lem:cc-reduction}. For example, assume $\mu_0$ is not a
 dominant integral weight. Let us show
\begin{equation}
    X_1[-n_1]\cdots X_l[-n_l]|\mu_0\rangle \tensor v \tensor v_\infty
    \equiv 0
    \mod{\gout^\orb},
\label{reduction:not-integrable:1}
\end{equation}
 where $X_i[-n_i] \in \hat\g^{(0)}$ ($X_i \in\g$, $n_i > 0$), $v\in M$,
 $v_\infty \in M^{(\infty)}_{\mu_\infty}$ and $\mod{\gout^\orb}$ denotes
 $\mod{\gout^\orb (M^{(0)}_{\mu_0} \tensor M \tensor
 M^{(\infty)}_{\mu_\infty})}$. (This abbreviation shall be used
 throughout this paper.)  Let $f_1(t)$ be an element of $\gout^\orb$
 such that $f_1(t) \sim X_1 \tensor t^{-n_1} + O(t^{n})$ for
 sufficiently large $n$. (Such $f_1$ exists due to the Riemann-Roch
 theorem. It is not difficult to construct such a function
 concretely.\footnote{A useful technique: for any $\g$-valued function
 $f(t)$ with poles in $\{0, Q_1,\dots,Q_L, \infty\}$, $f(t) +
 \Ad\gamma(f(\eps^{-1}t)) + (\Ad\gamma)^2 (f(\eps^{-2} t)) + \cdots +
 (\Ad\gamma)^{N-1} (f(\eps^{-N+1} t)) \in \gout^\orb$.}) Then we may
 replace $X_1[-n_1]$ by $\rho_0(f_1(t))$:
\begin{equation}
 \begin{split}
    &X_1[-n_1]\cdots X_l[-n_l]|\mu_0\rangle \tensor v \tensor v_\infty
\\
    =&
    \rho_0(f_1(t)) 
    X_2[-n_2]\cdots X_l[-n_l]|\mu_0\rangle \tensor v \tensor v_\infty
\\
    \equiv&
    - X_2[-n_2] \cdots X_l[-n_l]|\mu_0\rangle \tensor \\ & \tensor
    \left(\sum_{i=1}^L \rho_i(f_1(t)) v \tensor v_\infty
          + v \tensor \rho_\infty(f_1(t)) v_\infty\right)
    \mod{\gout^\orb}.
 \end{split}
\label{reduction:temp1}
\end{equation}
 By induction on $l$, the problem is reduced to showing 
\begin{equation}
    |\mu_0\rangle \tensor v \tensor v_\infty \equiv 0
    \mod{\gout^\orb},
\label{reduction:not-integrable:2}
\end{equation}
 which immediately follows from \lemref{lem:cc-reduction}. 

 To prove the second statement of \propref{prop:ccorb:integrable},
 assume that $\mu_0$ and $\mu_\infty$ are dominant integral
 weight. Then the irreducible quotients of $M^{(\ast)}_{\mu_\ast}$ are
 expressed as
\begin{equation}
    L^{(*)}_{\mu_*} =
    M^{(*)}_{\mu_*}/\sum_{i=0}^{N-1} 
    U(\n_-) f_i^{\langle \mu_*, \alpha_i^{\vee}\rangle + 1} |\mu_*\rangle.
\label{irred-quot}
\end{equation}
 See (10.4.6) of \cite{kac:90}. Therefore to prove
 \eqref{ccorb:integrable}, it is enough to show 
\begin{equation}
    U(\n_-) f_i^{\langle \mu_0, \alpha_i^{\vee}\rangle + 1} 
    |\mu_0\rangle \tensor M \tensor M^{(\infty)}_{\mu_\infty} \equiv 0
    \mod{\gout^\orb},
\end{equation}
 and a similar statement with the indices ``$0$'' and ``$\infty$'' for
 $\mu_*$, $M^{(*)}_{\mu_*}$ etc.\ interchanged. They are proved as above,
 namely by the arguments like \eqref{reduction:temp1} and
 \eqref{reduction:not-integrable:2}, because the weight of the singular
 vector $f_i^{\langle \mu_*, \alpha_i^{\vee}\rangle + 1} |\mu_*\rangle$
 is not a dominant integral weight.
\end{proof} 

\begin{proof}[Proof of \lemref{lem:cc-reduction}]
 Since $\kappa$ is not a dominant integral weights, there is an index
 $i$ ($0 \leqq i \leqq N-1$) such that $\kappa(\alpha_i^{\vee})$ is not
 a non-negative integer. By an easy calculation, we have
\begin{equation}
    e_i^n f_i^n v_\kappa = c v_\kappa, \qquad
    c 
    = n! \prod_{l=1}^n (\kappa(\alpha_i^{\vee}) - l + 1),
\label{eifivk}
\end{equation}
 for any $n \in \N$. Note that the constant $c$ never vanishes.

 Let $e(t)$ be an element of $\gout^\orb$ such that: (1) $e(t) \sim
 \phi_0^{-1}(e_i) + O(t^n)$ for sufficiently large $n$; (2)
 $\rho_\infty(e(t))v_\infty = 0$ (i.e., $e(t)$ has a zero of large order
 at $t=\infty$). Such an element can be constructed in the form $X
 \tensor F(t)$, where $X = E_{i,i+1}$ ($i=1,\dots,N-1$) or $X=E_{N,1}$
 ($i=0$) and $F(t)$ is a rational function. Hence we can rewrite
 $v_\kappa \tensor v \tensor v_\infty$ modulo $\gout^\orb$ as follows
 (cf.\ p.479 of \cite{tuy:89}):
\begin{equation}
 \begin{split}
    &v_\kappa \tensor v \tensor v_\infty
\\
    &=
    c^{-1} \rho_0(e_i^n f_i^n) v_\kappa \tensor v \tensor v_\infty
\\
    &=
    c^{-1} \rho_0(e(t)^n f_i^n) v_\kappa \tensor v \tensor v_\infty
\\
    &\equiv
    (-1)^n c^{-1} 
    \sum_{n_1 + \cdots + n_L = n}\frac{n!}{n_1! \cdots n_L!}
    \rho_0(f_i^n) v_\kappa \tensor
    \prod_{j=1}^L \rho_j(e(t))^{n_j} v \tensor v_\infty.
 \end{split}
\label{cc-reduction:temp}
\end{equation}
 Recall that $\rho_j(e(t)) = \rho_j(X \tensor F(t))$ is locally
 nilpotent on $M_j$ (Corollary 1.4.6 of \cite{tuy:89}). Thus the right
 hand side of \eqref{cc-reduction:temp} is $0$ for large $n$, which
 completes the proof of the lemma.
\end{proof}

Because of the difference of the sign in \eqref{def:phi0} and
\eqref{def:phi00} and the fact that $\Ad\beta$ is a Dynkin automorphism,
if $\tilde\lambda$ is a dominant integral weight for $\hat\g^{(0)}$,
$\tilde\lambda'$ is a dominant integral weight for $\hat\g^{(\infty)}$
and vice versa. (See \eqref{def:tilde-lambda}.)

\begin{corollary}
\label{cor:factor:integrable}
 Under the same assumption as \propref{prop:ccorb:integrable} we have
\begin{equation}
    \CC_\trig(M)
    \simeqq
    \bigoplus_{\lambda \in \wt(V), \tilde\lambda: \text{\rm dom. int.}}
    \CC_\orb(L^{(0)}_{\tilde\lambda} \tensor M 
     \tensor L^{(\infty)}_{\tilde\lambda'}),
\end{equation}
 where ``dom.\ int.'' means ``dominant integral weight''. Similarly the
 decomposition \eqref{factor:trig->orb} becomes
\begin{equation}
    \iota: \sCCtri(\M) \isoto
    \bigoplus_{\lambda \in \wt(V), \tilde\lambda: \text{\rm dom. int.}}
    \sCCorb(\L^{(0)}_{\tilde\lambda} \tensor \M 
    \tensor \L^{(\infty)}_{\tilde\lambda'}),
\label{factor:trig->orb:integrable}
\end{equation}
 where $\L^{(*)}_{\mu} = L^{(*)}_{\mu} \tensor \calO_{S_0}$ ($* =
 0,\infty$).
\end{corollary}

\section{Locally freeness}
\label{sec:loc-free}

The main theorem of this paper is proved in this section. We show that
$\sCC$ is locally free at the discriminant locus, $S_0 = \{q=0\} \subset
S$, provided that all modules inserted are {\em integrable highest
weight modules}. Thus, combining the result in \cite{kur-tak:97}, we
have locally freeness of the sheaf of conformal coinvariants and
consequently locally freeness of its dual, the sheaf of conformal
blocks. Corresponding statement for the Weyl modules has been proved in
\propref{prop:coherence}.

The arguments in this section is parallel to those in \S7.3 of
\cite{tsu-kuw:04}.

We assume the condition \eqref{M:finiteness} for $\M$. Moreover we
assume that all $M_i$'s are integrable highest weight modules. In
particular, the level $k$ is a non-negative integer.

\begin{theorem}
\label{thm:loc-free}
 The sheaf $\sCC(\M)$ and hence the sheaf $\sCB(\M)$ are locally free
 $\calO_{S}$-sheaves.
\end{theorem}

The rest of the paper is devoted to the proof of this theorem. The main
strategy of the proof is the same as that of \cite{tsu-kuw:04}. See
also \cite{shi-uen:99}, \cite{nag-tsu:02} and \cite{uen:95}:
\begin{itemize}
 \item Since $\sCC(\M)$ is coherent, it is sufficient to prove that each
       stalk $\sCC(\M)_s$ ($s \in S$) is a free $\calO_{S,s}$-module. 

 \item We proved the locally freeness of $\sCC(\M)$ on the non-singular
       part $S\setminus S_0$ in \cite{kur-tak:97}. Thus we have only to
       prove the case $s \in S_0$.

 \item When $s \in S_0$, we prove that the stalk of the completion of
       $\sCC(\M)$ at $s$ is isomorphic to $\sCCtri(\M)_s[[q]]$.

 \item $\sCC(\M)_s$ is a free $\calO_{S,s}$-module because of
       \propref{prop:loc-free:trig/orb} and the faithfully flatness of
       the completion functor.
\end{itemize} 

We define completion of the sheaf $\sCC(\M)$ along the divisor $S_0$ of
$S$ by taking completion of each ingredient of the definition
\eqref{def:sheaf-conf-coinvariant}. Let $\hatOS$ be the completion of
$\calO_S$ along $S_0$:
\begin{equation}
    \hatOS := \projlim_{n\to\infty} \calO_S/\m_{S_0}^n.
\label{def:hatOS}
\end{equation}
($\m_{S_0}$ is the definining ideal of $S_0$.) As an $\calO_{S_0}$-module,
it is isomorphic to the ring of formal power series: $\hatOS \simeqq
\calO_{S_0} [[q]]$.  The completion of $\M = M \tensor \calO_S$ is obviously
\begin{equation}
    \hatM := \M \tensor_{\calO_S} \hatOS \simeqq \M_{S_0} [[q]],
\label{hatM}
\end{equation}
where $\M_{S_0} := M \tensor \calO_{S_0}$. The $\calO_S$-Lie algebra $\goutX$
acts on $\hatM$ naturally as follows: a germ $f(P)$ of $\goutX$ at
$(q=0; Q_1,\dots,Q_L) \in S_0$ is expanded at
$(q=0;(0,0,0)_k;Q_1,\dots,Q_L) \in \tilde\X$ in terms of the coordinates
$(q;x_k,y_k)_k$ as
\begin{equation}
 \begin{split}
    f(P) &= \sum_{m,n = 0}^\infty f_{k,m,n}(s) x_k^m y_k^n
\\
         &= \sum_{m,n = 0}^\infty f_{k,m,n}(s) x_k^{m-n} q^n
          = \sum_{m,n = 0}^\infty f_{k,m,n}(s) y_k^{n-m} q^m
\\     
         &= \sum_{n = 0}^\infty f_{k,n,x}(s,x_k)  q^n
          = \sum_{m = 0}^\infty f_{k,m,y}(s,y_k)  q^m,
 \end{split}
\label{expansion:f:0,infty}
\end{equation}
where $f_{k,m,n}(s) \in \g\tensor\calO_{S_0}$  ($s\in S_0$) and
\begin{equation}
    f_{k,n,x}(s,x_k) = \sum_{m=0}^\infty f_{k,m,n}(s) x_k^{m-n}, \qquad
    f_{k,m,y}(s,y_k) = \sum_{n=0}^\infty f_{k,m,n}(s) y_k^{n-m}.
\label{def;fkm}
\end{equation}
Note that the periodicity condition \eqref{sec-goutX} implies
\begin{equation}
    \Ad\gamma (f_{k,m,n}(s)) = \eps^{m-n} f_{k,m,n}(s),\qquad
    \Ad\beta  (f_{k,m,n}(s)) = f_{k+1,m,n}(s),
\label{fkmn:period}
\end{equation}
and hence $f_{k,n,x}(s,x_k)$ and $f_{k,m,y}(s,y_k)$ are meromorphic
(rational) function on $\Proj^1$ with quasi-periodicity
\begin{equation}
 \begin{split}
     &f_{k,n,x}(s,\eps x_k) = \Ad\gamma (f_{k,n,x}(s,x_k)),
\\
     &f_{k,m,y}(s,\eps^{-1} y_k) = \Ad\gamma (f_{k,m,y}(s,y_k)),
\\
     &f_{k+1,n,x}(s,x_k) = \Ad\beta (f_{k,n,x}(s,x_k)),
\\
     &f_{k+1,m,y}(s,y_k) = \Ad\beta (f_{k,m,y}(s,y_k)),
 \end{split}
\label{fknx/y:period}
\end{equation}
the poles of which are in the divisor
\begin{equation}
    \ast \tilde D_0 
    := \sum_{i=1}^L\sum_{j\in\Z/N\Z} \ast [\eps^j Q_i]
     + \ast [0] + \ast [\infty].
\label{def:tildeD0}
\end{equation}
Therefore $\{f_{k,n,x}(s,x_k)\}_{k\in\Z/N\Z}$ and
$\{f_{k,n,y}(s,y_k)\}_{k\in\Z/N\Z}$, namely the $n$-th coefficients of
the expansion
\begin{equation}
    f(P) = \sum_{n=0}^\infty f_n(P) q^n,
\label{expansion:f}
\end{equation}
define a section $f_n(P)$ of $\gtwX|_{\pi_{\X/S}^{-1}(S_0)}$ with poles
at $Q_1,\dots,Q_L, 0, \infty$. See \eqref{def:gtwX} and
\eqref{sec-goutX}. The action of $f(P)$ on $v \tensor g \in \hatM$ ($v
\in \M$, $g \in \hatOS$) is defined by
\begin{equation}
    f(P)\cdot (v \tensor g)
    :=
    \sum_{i=1}^L \sum_{n=0}^\infty 
    \rho_i(f_n(P))v \tensor q^n g.
\label{def:f(P)(v*g)}
\end{equation}
Here $\rho_i$ denotes the usual action of the Laurent expansion of
$f_n(P)$ at $Q_i$ on $M_i$. The space of coinvariants of $\hatM$ with
respect to this action is the completion of $\sCC(\M)$:
\begin{equation}
    \hatCC(\M) := \hatM / \goutX (\hatM).
\label{def:hatCC}
\end{equation}

\begin{lemma}
\label{lem:completion}
\begin{equation}
    \hatCC(\M) \simeqq \sCC(\M)\tensor_{\calO_S}\hatOS.
\label{completion}
\end{equation}
\end{lemma}

\begin{proof}
 By definition, we have an exact sequence
\begin{equation}
    \goutX \tensor \M \to \M \to \sCC(\M) \to 0.
\label{goutXM->M->CC}
\end{equation}
 Tensoring $\hatOS$ we obtain an exact sequence
\begin{equation}
    (\goutX \tensor \M)[[q]] \to \hatM \to 
    \sCC(\M)\tensor\hatOS \to 0.
\label{completion:goutXM->M->CC}
\end{equation}
 It is sufficient to show that the image of the map $(\goutX \tensor
 \M)[[q]] \to \hatM$ is $\goutX (\hatM)$ defined by the
 action \eqref{def:f(P)(v*g)}. This is almost trivial since $f(P)\cdot v$
 for $f(P) \in \goutX$ and $v\in \M$ is expressed as
\begin{equation*}
    f(P) \cdot v = 
    \sum_{i=1}^L \sum_{n=0}^\infty 
    q^n \rho_i(f_n(P))v,
\end{equation*}
 because of the expansion \eqref{expansion:f}.
\end{proof}

Next step is to make a completion of the isomorphism
\eqref{factor:trig->orb:integrable}. For this purpose we need a lemma on
Verma modules of $\hat\g^{(0)}$ and $\hat\g^{(\infty)}$. Note that the
Verma module $M^{(*)}_\mu$ of $\hat\g^{(*)}$ ($*=0,\infty$) is graded
by the degree:
\begin{equation}
    M^{(*)}_\mu = \bigoplus_{d\geq 0} M^{(*)}_\mu(d), \qquad
    M^{(*)}_\mu(d) := \sum_{n_1+\cdots+n_l = d} 
    \C X_1[-n_1] \cdots X_l[-n_l] |\mu\rangle,
\label{Verma:degree}
\end{equation}
where $|\mu\rangle$'s are the highest weight vectors of $M^{(\ast)}_\mu$
and $X_i[-n_i] \in \hat\g^{(*)}$ ($X_i \in \g$, $n_i \in \N$).

\begin{lemma}
\label{lem:pairing-verma}
 (i)
 For any $\mu \in \h^*$ there exists a pairing between the Verma modules
 $M^{(0)}_{\mu\circ \Ad\beta}$ and $M^{(\infty)}_{-\mu}$:
\begin{equation}
    M^{(0)}_{\mu\circ\Ad\beta} \times M^{(\infty)}_{-\mu}: (u,v) \mapsto
    \left<u,v\right> \in \C,
\label{pairing}
\end{equation}
 which satisfies $\bigl< |\mu\circ\Ad\beta\rangle, |-\mu\rangle \bigr> =
 1$ and
\begin{equation}
    \left< X[n]u, v\right> + \left< u, \Ad\beta(X)[-n] v \right> = 0,
\label{pairing:invariance}
\end{equation}
 for any $u \in M^{(0)}_{\mu\circ\Ad\beta}$, $v \in
 M^{(\infty)}_{-\mu}$, $X\in\g$ and $n \in \Z$.

 (ii)
 $\left< M^{(0)}_{\mu\circ\Ad\beta}(n), M^{(\infty)}_{-\mu}(n')\right> =
  0$ if $n\neq n'$.

 (iii)
 The radical $R^{(0)} = \{ u \in M^{(0)}_{\mu\circ\Ad\beta} \mid \left<
 u, v \right> = 0 \text{\ for all } v \in M^{(\infty)}_{-\mu}\}$ is the
 largest proper submodule of $M^{(0)}_{\mu\circ\Ad\beta}$. Similarly the
 radical in $M^{(\infty)}_{-\mu}$ is the largest proper submodule. Hence
 the pairing descends to a non-degenerate pairing between the
 irreducible quotients $L^{(0)}_{\mu\circ\Ad\beta}$ and
 $L^{(\infty)}_{-\mu}$.
\end{lemma}

\begin{proof}
 Let $\nu$ be the anti-isomorphism
\begin{equation}
    \nu: U\hat\g^{(0)} \owns X[n] = X \tensor t^n
    \mapsto - X \tensor t^n = - X[-n] \in U\hat\g^{(\infty)},
    \qquad \nu(\hat k) = \hat k.
\label{def:nu:U->U}
\end{equation}
 This induces a linear isomorphism $\nu_\beta:
 M^{(0)}_{\mu\circ\Ad\beta} \to \Hom_\C(M^{(\infty)}_{-\mu},\C)$ defined
 by
\begin{equation}
    \nu_\beta(x|\mu\circ\Ad\beta\rangle) 
    = \langle -\mu| \nu(\Ad\beta(x)),
\label{def:nu-beta}
\end{equation}
 where $\langle -\mu|$ is the generating vector of the right
 $\hat\g^{(\infty)}$-module $\Hom_\C(M^{(\infty)}_{-\mu},\C)$,
 normalised by $\left< -\mu \mid -\mu \right> = 1$. We define the
 pairing $\left< , \right>$ by
\begin{equation}
    \left< v,v' \right> := \nu_\beta(v) v'.
\label{def:pairing}
\end{equation}
 Straightforward computation shows that for $x \in U\hat\g^{(0)}$ we have
\begin{equation}
    \left< x v , v' \right> = \left< v, \nu(\Ad\beta(x)) v' \right>,
\label{pairing:invariance:Ug}
\end{equation}
 which means \eqref{pairing:invariance} for $x = X[n]$.

 (ii) follows from the construction.

 (iii) is proved in the same way as Proposition 3.26 of \cite{wak:99}.

\end{proof}

Let $\{e_{\lambda,d,i}\}$ be a basis of $L^{(0)}_{\tilde\lambda}(d)$ and
$\{e_{\lambda,d}^{i}\}$ be its dual basis of
$L^{(\infty)}_{\tilde\lambda'}(d)$ with respect to $\left< , \right>$.

\begin{proposition}
\label{prop:hat-i}
 There exists an isomorphism
\begin{equation}
    \hat\iota:
    \hatCC(\M) \to 
    \bigoplus_{\lambda \in \wt(V), \tilde\lambda: \text{\rm dom. int.}}
    \sCCorb(\L_{\tilde\lambda} \tensor \M \tensor \L_{\tilde\lambda'})[[q]], 
\label{hat-i}
\end{equation}
 of $\hatOS$-modules defined by
\begin{equation}
    \hat\iota([v])
    :=
    \bigoplus_{\lambda}
    \left[\sum_{d=0}^\infty \sum_{i}
    e_{\lambda,d,i} \tensor v \tensor e_{\lambda,d}^i \right] q^d,
\label{def:hat-i}
\end{equation}
 for $v\in \M$.
\end{proposition}

\begin{proof}
 First we prove the well-definedness of \eqref{def:hat-i}, for which it
 is enough to show the well-definedness of its component $\iota_\lambda:
 \hatCC(\M) \to \sCCorb(\M_{\tilde\lambda} \tensor \M \tensor
 \M_{\tilde\lambda'})[[q]]$, namely, 
\begin{equation}
    \sum_{d=0}^\infty \sum_{i}
    e_{\lambda,d,i} \tensor f(P)\cdot v \tensor e_{\lambda,d}^i q^d
    \in
    \gout^\orb (\M_{\tilde\lambda} \tensor \M \tensor \M_{\tilde\lambda'})
\label{hat-i:well-def}
\end{equation}
 for $f(P) \in \goutX$ and $v \in \M$. Since
\begin{equation}
    f(P) \cdot
    \left(\sum_{d=0}^\infty \sum_{i}
    e_{\lambda,d,i} \tensor v \tensor e_{\lambda,d}^i 
    \right) q^d  \in
    \gout^\orb (\M_{\tilde\lambda} \tensor \M \tensor \M_{\tilde\lambda'}),
\label{hat-i:well-def:2}
\end{equation}
 we have only to show that the left hand side of \eqref{hat-i:well-def}
 is equal to the left hand side of \eqref{hat-i:well-def:2}, which is
 equivalent to an equation in $L_{\tilde\lambda} \tensor
 L_{\tilde\lambda'}$:
\begin{equation}
    \sum_{d=0}^\infty \sum_{i}\left(
    \rho_0(f(P)) \cdot e_{\lambda,d,i} \tensor e_{\lambda,d}^i 
    +
    e_{\lambda,d,i} \tensor \rho_\infty(f(P)) \cdot e_{\lambda,d}^i 
    \right) q^d
    = 0.
\label{hat-i:invariance}
\end{equation}
 Recall that the germ of $f(P)$ at $0$ and the germ at $\infty$ is
 related by $f(P)_\infty = \Ad\beta (f(P)_0)$. See \eqref{sec-goutX} of
 this paper or (18) of \cite{tak:04}. Hence using the expansion
 \eqref{expansion:f:0,infty} and the invariance \eqref{pairing:invariance} of
 the pairing, we can show \eqref{hat-i:invariance} in the same way as
 the proof of Claim 3 in the proof of Theorem 6.2.1 in
 \cite{tuy:89}. Thus $\hat\iota$ is well-defined.

 Obviously the $q=0$ part of $\hat\iota$ is the isomorphism $\iota$,
 \eqref{factor:trig->orb:integrable}. Therefore by termwise
 approximation (in analytic language) or, in other words, by Nakayama's
 lemma (in algebraic language), $\hat\iota$ is shown to be an
 isomorphism. 
\end{proof}

With these preparations, the proof of the locally freeness of $\sCC(\M)$
goes as follows. As is mentioned after the statement of
\thmref{thm:loc-free}, it is enough to prove that the stalk $\sCC(\M)_s$
at $s \in S_0$ is a free $\calO_{S,s}$-module. Since $\hatCC(\M)_s$ is
isomorphic to $\sCC(\M)_s \tensor_{\calO_{S,s}} \tensor \hatOS{}_{,s}$
(\lemref{lem:completion}) and to $\bigoplus_\lambda
\sCCorb(\L_{\tilde\lambda} \tensor \M \tensor \L_{\tilde\lambda'})_s
[[q]]$ 
(\propref{prop:hat-i}), we have an isomorphism
\begin{equation}
    \sCC(\M)_s \tensor_{\calO_{S,s}} \hatOS{}_{,s}
    \simeqq
    \bigoplus_\lambda
    \sCCorb(\L_{\tilde\lambda} \tensor \M \tensor \L_{\tilde\lambda'})_s
    [[q]].
\label{cc=sum-ccorb:completed}
\end{equation}
The right hand side of \eqref{cc=sum-ccorb:completed} being a free
$\hatOS{}_{,s}$-module (\propref{prop:loc-free:trig/orb} (ii) and
\propref{prop:ccorb:integrable}), faithfully flatness of $\hatOS{}_{,s}$
over $\calO_{S,s}$ implies that $\sCC(\M)_s$ is a free
$\calO_{S,s}$-module. Thus \thmref{thm:loc-free} is proved.

\section{Concluding comments}
\label{sec:conclusion}

We have proved locally freeness of $\sCC(\M)$ in two cases; Weyl module
case (\propref{prop:coherence}) and integrable highest weight module
case (\thmref{thm:loc-free}). A few comments are in order:

\begin{itemize}
 \item In the Weyl module case, $\sCC(\M) \simeqq V \tensor \calO_S$ as
       shown in the proof of \propref{prop:coherence} and the rank of
       $\sCC(\M)$ is $\dim V$.

 \item In the integrable highest module case, the rank is computed by
       further degenerating the orbifold. Degeneration of the type
       $Q_i \to Q_j$ should be considered in the same way as in
       \cite{tuy:89} or \cite{nag-tsu:02}. The final results of the
       degeneration is a combination of the three-punctured orbifold
       $\Proj^1/C_N$. In principle a Verlinde-type formula would be
       obtained in this way.

 \item In \cite{kur-tak:97} we have shown that $\sCC(\M)$ has a flat
       connection. It has a regular singularity along $S_0 = \{q=0\}$,
       which is easily deduced from the explicit form of the connection,
       Theorem 5.9 of \cite{kur-tak:97}. Hence there is a one-to-one
       correspondence between flat sections around $S_0$ and its
       restriction to $S_0$ or, in other words, the ``initial value'' at
       $S_0$ because of the locally flatness.
\end{itemize}

\section*{Acknowledgments}

The author expresses his gratitude to Akihiro Tsuchiya who explained
details of \cite{tuy:89} and \cite{nag-tsu:02}, Toshiro Kuwabara who
showed the manuscript of \cite{tsu-kuw:04} (the best guide to
\cite{tuy:89} for $\Proj^1$ case) before publishing, Michio Jimbo, Gen
Kuroki, Tetsuji Miwa, Hiroyuki Ochiai, Kiyoshi Ohba, Nobuyoshi Takahashi,
Tomohide Terasoma and Yasuhiko Yamada for discussion and comments.

The atmosphere and environment of Institute for Theoretical and
Experimental Physics (Moscow, Russia) and the conference ``Infinite
Dimensional Algebras and Integrable Systems'' (Faro, Portugal) were very
important. The author thanks their hospitality.


\begin{thebibliography}{TUY}
\bibitem[BD]{bel-dri:82}
A. A. Belavin, V. G. Drinfeld,
Solutions of the classical Yang-Baxter equations for simple Lie algebras.
{\it Funkts. Anal. i ego Prilozh.} {\bf 16-3},
1--29
(1982)
(in Russian);
{\it Funct. Anal. Appl.} {\bf 16},
159--180
(1982)
(English transl.)

\bibitem[E]{eti:94}
P. I. Etingof,
Representations of affine Lie algebras, elliptic $r$-matrix
systems, and special functions.
{\it Comm. Math. Phys.} {\bf 159},
471--502
(1994).

\bibitem[FFR]{f-f-r:94}
B. Feigin, E. Frenkel, N. Reshetikhin,
Gaudin model, Bethe Ansatz and critical level.
Commun. Math. Phys. {\bf 166},
27--62
(1994)

\bibitem[K]{kac:90}
V. G. Kac,
{\it Infinite dimensional Lie algebras},
3rd Edition,
Cambridge University Press 1990.

\bibitem[KL]{kaz-lus:93}
D. Kazhdan, G. Lusztig,
Tensor structures arising from affine Lie algebras. I, II,
{\it J. Amer. Math. Soc.} {\bf 6},
905--948, 949--1011
(1993).

\bibitem[KT]{kur-tak:97}
G. Kuroki, T. Takebe,
Twisted Wess-Zumino-Witten models on elliptic curves.
{\it Comm. Math. Phys.} {\bf 190},
1--56
(1997).

\bibitem[NT]{nag-tsu:02}
K. Nagatomo, A. Tsuchiya,
Conformal field theories associated to regular chiral vertex operator
algebras I: theories over the projective line,
{\tt math.QA/0206223}.

\bibitem[SU]{shi-uen:99}
Y. Shimizu, K. Ueno,
{\it Moduli theory III},
(Iwanami, Tokyo, 1999)
Gendai Suugaku no Tenkai series
(in Japanese);
{\it Advances in moduli theory},
Translations of Mathematical Monographs, {\bf 206},
Iwanami Series in Modern Mathematics,
American Mathematical Society, Providence, U.S.A. (2002)
(English translation)

\bibitem[T]{tak:04}
T. Takebe,
Trigonometric Degeneration and Orbifold Wess-Zumino-Witten Model. I 
In {\it the Proceedings of the 6th International workshop on Coformal
and Integrable models, Chernogolovka, Sep. 2002},
{\it International Journal of Modern Physics, A}, {\bf 19}, 
Supplement,
418--435
(2004)

\bibitem[TK]{tsu-kuw:04}
A. Tsuchiya, T. Kuwabara,
Introduction to Conformal Field Theory,
to appear as MSJ Suugaku Memoir of Mathematical Society of Japan.

\bibitem[TUY]{tuy:89}
A. Tsuchiya, K. Ueno, Y. Yamada,
Conformal field theory on universal family of stable curves with
gauge symmetries.
In
{\it Integrable systems in quantum field theory and statistical mechanics,
Adv. Stud. Pure Math.} {\bf 19},
459--566
(1989).

\bibitem[U]{uen:95}
K. Ueno,
On conformal field theory,
In {\it Vector bundles in algebraic geometry (Durham, 1993)},
ed. by N. J. Hitchin, P. E. Newstead and W. M. Oxbury,
{\it London Math. Soc. Lecture Note Ser.} {\bf 208}, 
(Cambridge Univ. Press, Cambridge, 1995)
pp. 283--345,

\bibitem[Wa]{wak:99}
M. Wakimoto,
{\it Infinite-dimensional Lie algebras},
(Iwanami, Tokyo, 1999)
Gendai Suugaku no Tenkai series
(in Japanese); 
Translations of Mathematical Monographs, {\bf 195},
Iwanami Series in Modern Mathematics,
American Mathematical Society, Providence, U.S.A. (2001)
(English translation by K. Iohara)


\bibitem[Wo]{wol:83}
S. Wolpert,
On the homology of the moduli space of stable curves.
{\it Ann. of Math.} {\bf 118},
491--523
(1983).

\bibitem[Z]{zam:89}
A. B. Zamolodchikov,
Exact solutions of conformal field theory in two dimensions and critical
phenomena. 
{\it Rev. Math. Phys.} {\bf 1}
197--234
(1989).
(Translated from the Russian by Y. Kanie.)

\end{thebibliography}
\end{document}